\documentclass[letterpaper,10pt]{amsart}
 \DeclareFontFamily{U}{wncy}{}
    \DeclareFontShape{U}{wncy}{m}{n}{<->wncyr10}{}
    \DeclareSymbolFont{mcy}{U}{wncy}{m}{n}
    \DeclareMathSymbol{\Sh}{\mathord}{mcy}{"58}
\usepackage{scalerel}
\usepackage[mathscr]{euscript}

\usepackage{hyperref}
\usepackage{amssymb}
\usepackage{amscd}
\usepackage{amsxtra}
\usepackage{mathrsfs}
\usepackage{tikz}
\usepackage{tikz-cd}
\usetikzlibrary{matrix}
\newcommand{\rank}{\mathrm{rank}\;}
\newcommand{\set}[1]{\left\{ #1\right\}}

\newcommand{\CC}{\mathbb{C}}

\newcommand{\FF}{\mathbb{F}}

\newcommand{\NN}{\mathbb{N}}
\newcommand{\PP}{\mathbb{P}}

\newcommand{\ZZ}{\mathbb{Z}}
\newcommand{\ip}[1]{\left\langle #1 \right\rangle }

\newcommand{\ep}{\varepsilon}

\newcommand{\si}{\sigma}

\newcommand{\Gg}{\mathcal{G}}
\newcommand{\Hh}{\mathcal{H}}

\newcommand{\Ll}{\mathcal{L}}

\newcommand{\Oo}{\mathcal{O}}

\usepackage{bbm}

\newtheorem{thm}{Theorem}[section]
\newtheorem{df}[thm]{Definition}

\newtheorem{lem}[thm]{Lemma}
\newtheorem{prp}[thm]{Proposition}

\begin{document}
\title{Weil-Chatelet Groups of Rational Elliptic Surfaces}
\author{Nadir Hajouji}
\maketitle

\begin{abstract}
    We classify pairs $(S, \gamma)$, consisting of a rational elliptic surface $S$ and a Galois cover $\gamma$ of the base, which satisfy a condition we call $\Ll$-stability. 
    We explain how to use the theory of Mordell-Weil lattices to compute the kernel of the  restriction maps of Weil-Chatelet groups for $\Ll$-stable pairs.
    We also prove results about the injectivity of restriction maps of Weil-Chatelet groups for some pairs which are not $\Ll$-stable.
\end{abstract}

\section{Introduction}

Let $K$ be a perfect field, $\overline{K}$ the algebraic closure and $\Gg = Gal(\overline{K}/K)$.
Let $E/K$ be an elliptic curve.
The Weil-Chatelet group $E/K$ is the Galois cohomology group:
\[ WC(E/K) := H^1(\Gg, E(\overline{K})) \]

There are various reasons one may be interested in the Weil-Chatelet group of an elliptic curve:
\begin{itemize}
    \item The Weil-Chatelet group is usually introduced\footnote{At least that is how it is introduced in \cite{sil1}.} as an ``obstruction" to computing the Mordell-Weil group.
    \item The Weil-Chatelet group is isomorphic to a quotient of the Brauer group of the elliptic curve, and the Weil-Chatelet group contains the Tate-Shafarevich group of the elliptic curve. Sometimes one is interested in $WC$ because of its relationship to these groups. 
    \item However, the group is interesting in it's own right - classes in the Weil-Chatelet group represent isomorphism classes of genus one curves without a rational point, and those genus one curves can be hard to study ``directly". 
\end{itemize}
This work is motivated by the latter - we want to understand $WC(E/K)$ because we want to be able to classify genus one fibrations without a section.

Unfortunately, the Weil-Chatelet group is \emph{notoriously} difficult to compute.
\footnote{In most cases, it is impossible to even describe the whole group, since it is not finitely generated.} 
In fact, the mere problem of computing a kernel $WC(E/K) \to WC(E/L)$,
where $L/K$ is a finite Galois extension,
can be hard to do in general, because it requires us to understand $E(L)$ as a $Gal(L/K)$-module, which again requires us to understand the Mordell-Weil group,
and there is no general algorithm for computing generators of the Mordell-Weil group.

However, there is one setting where Mordell-Weil groups can be computed without much trouble: when $E/K$ is the generic fiber of a rational elliptic surfaces.

\begin{itemize}
    \item Right off the bat, we know that the Mordell-Weil group has rank at most 8, and in fact the torsion-free part of the Mordell-Weil group can be embedded into the $E_8$ lattice.
    
    \item We can determine the precise rank of $E/K$ using the Shioda-Tate formula as soon as we know the isomorphism types\footnote{In fact, we only need to know the number of irreducible components in each fiber.} of the singular fibers in a minimal, smooth model of $E/K$.
    
    \item In fact, it is easy to find an explicit generating set for the Mordell-Weil group of a rational elliptic surface - the group is generated by the sections of height at most 2, and one can classify all sections of height at most 2 using a computer algebra system. 
\end{itemize}

The goal of this work is to exploit these ideas to better understand Weil-Chatelet groups of certain elliptic surfaces.

\subsection{Summary of Main Results}

Let $\pi:S\to \PP^1$ be a smooth, minimal, non-isotrivial elliptic surface,
let $E/K$ be the generic fiber,
and let $\Ll_{S/\PP^1}$ be the fundamental line bundle of the elliptic surface.
Let $\gamma : \PP^1 \to \PP^1$ be a Galois cover, and let $S_\gamma \to \PP^1$ be the minimal, smooth resolution of the base changed elliptic surface $S \times_\gamma \PP^1 \to \PP^1$.

{\df{We say that the pair $(S, \gamma)$ is $\Ll$-stable if $\Ll_{S/\PP^1} \cong \Ll_{S_\gamma/\PP^1}$.
The index of an $\Ll$-stable pair is the degree of the Galois cover $\gamma$.}}

If we have an $\Ll$-stable pair $(S,\gamma)$, 
where $S$ is rational,
then the base-changed surface is also rational.
Thus, we can use the theory of Mordell-Weil lattices of rational elliptic surfaces to describe the Mordell-Weil group of $E$ before \emph{and} after base change.
This means can determine the isomorphism type of $E$ after base change not only as an abelian group, but also as a Galois module. 

{\df{Let $(S,\gamma)$ be an $\Ll$-stable pair.
\begin{itemize}
    \item 
We say that $(S,\gamma)$ is \emph{rank-stable} if the rank of the Mordell-Weil group does not increase after base-changing.
\item 
We say that $(S,\gamma)$ is $WC$-stable if the induced map $WC(E/K)\to WC(E/K)$ is injective.

\end{itemize}
}}

We will show in Theorem \ref{thm:reptheory} that rank-stable implies $WC$-stable\footnote{At least for the pairs we are studying.}.
Rank stability is easier to understand than $WC$-stability, and in many cases,
it is enough to determine whether we have $WC$-stability.
However, we will see that the two conditions are not equivalent - we will see two examples of $\Ll$-stable pairs which are not rank stable, but which are $WC$-stable.

\begin{itemize}
    \item We will classify all $\Ll$-stable pairs $(S,\gamma)$ where $E$ is the generic fiber of a rational elliptic surface, and $\gamma$ is a Galois cover of index $p \geq 5$.
    The assumption $p \geq 5$ is significant - there would be a lot more examples of $\Ll$-stable pairs if we allowed $p = 2,3$.
    However, our goal is to ultimately understand genus one fibrations without a section, and there are more geometric methods we can use to understand those fibrations.
    
    \item Once we've classified $\Ll$-stable pairs, we use Shioda-Tate to determine which pairs are rank stable.
    We will show that those pairs are also $WC$-stable, so we can essentially ignore them.
    That will leave only a handful of $\Ll$-stable pairs which could be hiding a torsor.

    \item Finally, we will determine which of the rank-unstable pairs are $WC$-stable.
    \begin{itemize}
        \item For some of the pairs $(S,\gamma)$, $E$ is the generic fiber of an \emph{extremal} rational elliptic surface.
        For these pairs, the map $WC(E/K) \to WC(E/K)$ is not injective, and in fact we can show that the kernel is either isomorphic to $\ZZ/5$ or $\ZZ/5 \times \ZZ/5$.
        \item There are only two pairs $(S,\gamma)$ where the original elliptic surface is not extremal, but the pair is not rank stable. 
        For these pairs, we will actually compute the Mordell-Weil groups before and after base change explicitly and use that to prove $WC$-stability.

    \end{itemize}
\end{itemize}

As summary of these results can be found in Table \ref{tab:summary}.

\subsection{Structure of the paper}

The paper is structured as follows:
\begin{itemize}
    \item In Section \ref{sec:galreps}, we discuss $WC(E/K)$ abstractly and explain how we intend to compute the kernel of $WC(E/K)\to WC(E/K)$ using data about the Mordell-Weil group after base change. 
    \item In Section \ref{sec:numconstraints}, we translate the $\Ll$-stability condition into a constraint on the fibers of the original elliptic surface $S\to \PP^1$. 
    We classify all $\Ll$-stable pairs in this section.
    \item In Section \ref{sec:examples}, we compute the kernel of $WC(E/K) \to WC(E/K)$ for all $\Ll$-stable pairs we found in the previous section.
    
\end{itemize}

\section{Galois Representations and $WC$}
\label{sec:galreps}

\subsubsection{Notation}
In this section, we use the following notation:

\begin{itemize}
    \item $K$ is a field of characteristic 0.\footnote{In later sections, we will be working with $K = k(t)$,
    where $k$ is an algebraically closed field of characteristic 0.}
    \item $\overline{K}$ is the algebraic closure of $K$.
    \item $\Gg = Gal(\overline{K}/K)$ is the absolute Galois group.
    \item $K'/K$ is a finite, Galois extension, with $Gal(\overline{K}/K') = \Hh$.
    Note that $\Hh$ is a normal subgroup of $\Gg$.
    \item $G = Gal(K'/K)$.
    \item  $E/K$ is an elliptic curve.
    We assume that $E$ is given by a fixed Weierstrass equation,
so that points in $E(K')$ (resp. $E(\overline{K})$) are endowed with the structure of a $G$-module (resp. $\Gg$-module).
\end{itemize}

\subsection{Galois Cohomology}

Let $WC(E/K) = H^1(\Gg,E(\overline{K}))$ and $WC(E/K') = H^1(\Hh, E(\overline{K}))$.
Then the Hochschild-Serre spectral sequence in group cohomology (see e.g. Cor.6.7.4 in \cite{poonenqpoints}) implies that the following sequence is exact:
\[ 0 \to H^1(G,E(K')) \to WC(E/K) \to WC(E/K')\]
Thus, we may identify $H^1(G,E(K'))$ with the kernel of $WC(E/K)\to WC(E/K')$.
This reduces the problem of determining $WC$-stability to understanding the group $H^1(G,E(K'))$.

To compute $H^1(G, E(K'))$, we use the description as cocycles modulo coboundaries\footnote{See Appendix $B$ in \cite{sil1}, e.g.}.
A cocycle $\xi \in Z^1(G, E(K'))$ is a function $\xi : G \to E(K')$ that satisfies:
\[ \xi(\si_1 \si_2 ) = \xi(\si_1) + \si_1 \cdot \xi(\si_2) \qquad (\forall \si_1, \si_2 \in G) \]

A coboundary $\xi_P \in Z^1(G,E(K'))$ is a cocycle $\xi_P : G \to E(K')$ of the form:
\[ \xi_P(\si) = \si\cdot P-P \]
for some fixed $P \in E(K')$.

The group of coboundaries is relatively straightforward to understand: it is isomorphic, as a $G$-module,
to $E(K')/E(K)$.\footnote{There is an obvious surjection $E(K') \to B^1(G,E(K'))$, and the kernel of that surjection is $E(K')^G = E(K)$.}
The group of cocycles can be a little bit more complicated to deal with in general.
However, if we assume that $G$ is cyclic, say $G = \ip{\si}$, then any cocycle $\xi \in Z^1(G,E(K'))$ is determined by $\xi(\si)$:
the cocycle condition forces $\xi(\si^2) = \si \cdot \xi(\si) + \xi(\si)$,
$\xi(\si^3) = \si \cdot \xi(\si^2)+\xi(\si)$, etc.
In order for this to be well-defined, we need $\xi(\si^n) = 0$.
This turns out to be equivalent to requiring $\xi(\si)$ to be in the kernel of the trace map $E(K')\to E(K)$.
Thus, we can identify $Z^1(G,E(K'))$ with $\ker(E(K')\to E(K))$.

If $P \in \ker(E(K') \to E(K))$, then the cocycle determined by $\xi(\si) = P$ is a coboundary if and only if there exists $Q \in E(K')$ such that $P = \si\cdot Q - Q$.
\begin{equation}
\label{eq:kerquotient}
H^1(G,E(K')) = Z^1(G,E(K'))/B^1(G,E(K')) \cong \ker(Tr:E(K')\to E(K))/\mathrm{Im}\;(Q\mapsto (1-\si))Q) \end{equation}

We will compute the kernel of $WC(E/K) \to WC(E/K')$ in two examples by directly computing the right hand side of \ref{eq:kerquotient} in the last two examples Section \ref{sec:examples}. 
This will require finding a generating set for $E(K')$,
which we can do using Theorem 8.33 in \cite{shiodaschuttmordellweillatticesbook}.

We can obtain more powerful results that do not require knowledge of a generating set of $E(K')$ if:
\begin{itemize}
    \item The Galois group $G$ is cyclic of prime order.
    \item The Mordell-Weil group after base changing is torsion-free.
\end{itemize}
We discuss this next.

\subsection{Cyclic Galois Groups of Prime Order}

When $G$ is cyclic \emph{of prime order}, things are even better, since we can appeal to the results in \cite{reinerintegralrepspgroups} to determine
$H^1(G,E(K'))$ up to isomorphism without actually needing to find a generating set for $E(K')$.

Precisely:

{\thm{\label{thm:reptheory} Let $E/K$ be an elliptic curve and $K'/K$ a cyclic extension of index $p$.
Assume that $E(K')$ is finitely generated and torsion-free.
Then:
\begin{enumerate}
    \item $E(K') = E(K) \oplus M$ as a $\ZZ[G]$-module.
    \item The rank of $M$ as an abelian group is divisible by $p-1$.
    \item If the rank of $E(K)$ is equal to the rank of $E(K')$, then:
    \begin{itemize}
        \item $E(K) = E(K')$.
        \item The trace map $E(K')\to E(K)$ is injective.
        \item The restriction map $WC(E/K)\to WC(E/K')$ is injective.
    \end{itemize}
\end{enumerate}
}}

This theorem follows directly from the results in  \cite{reinerintegralrepspgroups}.
We will need the following notation:
\begin{itemize}
    \item Let $G = \ZZ/p\ZZ$ and let $\ZZ[G]$ be the integral group ring.
    \item Let $\tau = \sum_{\si^i \in G} \si^i \in \ZZ[G]$ be the trace operator (that is, $\tau \cdot Q = Tr(Q)$ for any $Q \in E(K')$.)
    \item Let $\zeta_p = e^{2\pi i /p}$ and $\Oo = \ZZ[\zeta_p]$.
    Note that $\Oo$ is a $\ZZ[G]$-module\footnote{If we fix a generator $\si$ of $G$,
    then $g$ acts on an element of $\Oo$ by multiplication-by-$\zeta_p$.}.
    Furthermore, $\ZZ[G]/\ip{\tau} \cong \Oo$ as a $\ZZ[G]$-module.
    
\end{itemize}

\begin{proof}
The first point follows from the main result of \cite{reinerintegralrepspgroups}.

To prove the second point, consider the map $M\to M$ given by $m\mapsto (1-\si)M$.
Since $M^G = 0$, this map is injective, so $(1-\si)M$ is a submodule of $M$ which has the same rank as $M$.
Since every element of $(1-\si)M$ is annihilated by $\tau = 1+\si + \cdots + \si^{p-1}$, the $\ZZ[G]$-module structure of $(1-\si)M$ descends to a $\ZZ[G]/\ip{\tau}$-module structure.

But $\ZZ[G]/\ip{\tau} \cong \ZZ[ \zeta_p]$, where $\zeta_p$ is a primitive $p$th root of unity.
Since $(1-\si)M$ is finitely generated and torsion-free, it is isomorphic to a direct sum of fractional ideals of $ \ZZ[ \zeta_p]$, and these all have rank $p-1$ as abelian groups.

If the rank of $E(K)$ is equal to the rank of $E(K')$,
then $M = 0$ so $E(K') = E(K)$ by (1).
Thus, the trace map $E(K')\to E(K)$ is just the multiplication-by-$p$
map on $E(K)$.
Since $E(K)$ is torsion-free, every multiplication map is injective,
so $Tr:E(K')\to E(K)$ is injective.
Thus, $Z^1(Gal(K'/K), E(K')) = 0$, so $WC(E/K)\to WC(E/K')$ is injective.
\end{proof}

We can also use the results in \cite{reinerintegralrepspgroups} to prove that $WC(E/K)\to WC(E/K')$ in some cases:

{\thm{\label{thm:ext2nonext} Let $E/K$ be an elliptic curve with $E(K) = 0$,
and let $K'/K$ be a Galois extension of index $p \in \set{5,7,11,13,17,19,23}$. 
If $\mathrm{rank}\; E(K') \neq 0$, then $WC(E/K)\to WC(E/K')$ is not injective.

}}

In fact, we will show that one can determine the dimension of $H^1(G,E(K'))$ as an $\FF_p$-vector space
directly from the rank of $E(K')$ as an abelian group.

\begin{proof}
First, note that since $E(K) = 0$,
$Tr: E(K')\to E(K)$ is the zero map,
so the kernel of the trace map is all of $E(K')$.
Note that this means the $\ZZ[G]$-module structure of $E(K')$ descends to an $\Oo$-module structure.

Let $M$ denote the torsion-free part of $E(K')$.
Then $M$ is a regular $\Oo$-module.
Furthermore, since $p \leq 23$, $\Oo$ is a PID,
so $M \cong \Oo^r$ for some $r$.
Note that $r>0$ since we are assuming the pair is not rank stable.

Now, the kernel of the trace map is all of $E(K')$,
and the group of coboundaries is $(1-\si)E(K')$.
The action of $(1-\si)$ on $E(K')$ is the same as the action of multiplication by $(1-\zeta_p)$ on $\Oo^r$,
so the quotient $E(K')/(1-\si)E(K')$ is isomorphic to:
\[H^1(G,E(K')) \cong  \Oo^r/(1-\zeta_p) \Oo^r \cong (\Oo/(1-\zeta_p) \Oo)^r \cong (\ZZ/p\ZZ)^r \]

Since we are assuming the pair is not rank stable, it's clear that $r \neq 0$,
so $WC(E/K) \to WC(E/K')$ is not injective,
and the rank of the kernel as an $\FF_p$ vector space is equal to the rank of $M$ as an $\Oo$-module.

\end{proof}

{\thm{\label{thm:nonextremalwcstable}

Let $K$ be a field, $K'/K$ a cyclic Galois extension with Galois group $G = \ip{\si}$ (where $|G| = n$, say), $E/K$ be an elliptic curve and assume:
\begin{itemize}
\item{$E(K) \cong \ZZ$ as an abelian group.}
\item{$E(K')$ is generated, as a $\ZZ[G]$-module, by a single point $Q$.}
\item $Tr(Q)$ is a generator of $E(K)$.

\end{itemize}
Then $WC(E/K) \to WC(E/K')$ is injective.

}}

\begin{proof}

Since $Q$ generates $E(K')$ as a $\ZZ[G]$ module,
every element of $E(K')$ has the form:
\[  \sum_{\si^\ell \in G}  a_\ell (\si^\ell \cdot Q) \qquad (a_\ell \in \ZZ) \]
 
Since $Tr(\si^\ell Q) = Tr(Q) = P$, we compute:
\[ Tr\left( \sum_{\si^\ell \in G} a_\ell (\si^\ell \cdot Q) \right) = \sum a_\ell Tr(\si^\ell \cdot Q) = \sum a_\ell Tr(Q) =\left( \sum a_\ell \right)\cdot P\]
We know that $P$ is not a torsion point., so $ \sum_{\si^\ell \in G}  a_\ell (\si^\ell \cdot Q)$ has vanishing trace if and only if $\sum a_\ell = 0$.

Now, fix a point $Q_0 =  \sum_{\si^\ell \in G}  a_\ell (\si^\ell \cdot Q)$ in the kernel of the trace map and define:
\[ b_0 = 0 \quad b_{\ell+1} = a_\ell  +b_\ell \]
Let:
\[ Q_1 = \sum_{\si^\ell \in G} b_\ell (\si^\ell \cdot Q) \]
Then:
\[ (\si - 1) \cdot Q_1=  \sum_{\si^\ell \in G}(b_{\ell+1}- b_\ell)  (\si^\ell \cdot Q) = \sum_{\si^\ell \in G} a_\ell (\si^\ell \cdot Q) = Q_0 \]

This shows that every point $Q_0$ in the kernel of the trace map has the form $(1-\si) \cdot Q_1$, which proves $WC(E/K)\to WC(E/K')$ is injective.

\end{proof}

\section{$\Ll$-stable pairs over $\PP^1$}
\label{sec:numconstraints}

In this section, we classify $\Ll$-stable pairs over $\PP^1$ of index $p \geq 5$.
For an elliptic surface over $\PP^1$ and a Galois cover $\gamma: \PP^1\to \PP^1$,
it is easy to decide whether we have $\Ll$-stability if we know the fiber configuration before base change,
and we know where the Galois extension $\PP^1\to\PP^1$ is ramified - we can easily deduce the fiber configuration of the minimal resolution \emph{after} base changing.
 Once we know the fiber configuration after base change, we can determine the fundamental line bundle and decide whether the pair is $\Ll$-stable.
 
 If the minimal resolution after base changing is rational, we can also use the fiber configuration to determine the exact rank of the Mordell-Weil group using the Shioda-Tate formula (and exploiting the fact that $\rho(S) = 10$ for rational elliptic surfaces).
    This will allow us to determine when we have rank stability purely in terms of the fiber configuration and the the location of the ramification points in the next section.

Before discussing singular fibers of $\pi : S\to \PP^1$, we show that the Galois cover $\gamma: \PP^1\to \PP^1$ must be ramified over exactly two points,
as we will be using that fact repeatedly.
We will use the letter $C$ to refer to the domain of $\gamma$ and $\PP^1$ to refer to the codomain, to make it clear which genus zero curve we are talking about. 

\subsection{Constraints on $\gamma$}

Let $\gamma : C \to \PP^1$ be a Galois cover of degree $p$.
We want $C$ to have genus zero if we want to use this cover to construct $\Ll$-stable pairs.
The Riemann-Hurwitz formula (e.g. IV,Cor. 2.4 of \cite{hartshorne}) forces:
\[ 2(g(C) - 1) = 2p(g(\PP^1)-1)+ \deg R \]
where $\deg R$ is the degree of the ramification divisor.
We want $g(C) = 0$, so we solve for the degree of $R$:
\[-2 = -2p + \deg R \quad \implies\quad  2p-2 = \deg R \]
Since we have tame ramification, $\deg R = (p - 1)$ times the number of ramification points.
We want $g(C) = 0$, so we need $\deg R = 2p-2$, so there should be exactly two ramification points.

\subsection{General Constraints on Fiber Configurations}

We use coordinates $[t_0:t_1]$ on $\PP^1$, and refer to the points $[0:1],[1:0]$ as $0,\infty$.
We assume that the ramification points of $\gamma : C \to \PP^1$ are at $0,\infty$.

Let $\nu_0,\nu_\infty$ be the valuations at $0,\infty$ on $\PP^1$,
and $\nu_0',\nu_\infty'$ the valuations at the points lying over $0,\infty$ on $C$.

Let $\pi: S \to \PP^1$ be an elliptic surface with $\Ll_{S/\PP^1} \cong \Oo(d)$ for some positive integer $d$.
There is a Weierstrass model associated to $S\to \PP^1$ of the form:
\[ y^2 = x^3 + f_{4d} (t_0,t_1) x + g_{6d} (t_0,t_1)\]
where $f,g$ are homogenous polynomials whose degree is indicated by the subscript.

Let $\Delta = -16(4f^3+27g^2)$ be the discriminant of $S\to \PP^1$.
This is a homogenous polynomial of degree $12d$.
Let $\ep = 12d - \nu_0(\Delta)+\nu_\infty(\Delta)$.

If $\gamma: C\to \PP^1$ is a Galois cover ramified over $0,\infty$,
the Weierstrass equation after base changing may not be minimal (and in fact it will never be in the cases we care about).
Let $\tilde{f}, \tilde{g}$ be the Weierstrass coefficients of the minimal integral model after base changing,
and $\tilde{\Delta} = -16(4\tilde{f}^3+27\tilde{g}^2)$.
This allows us to compute the degree of the fundamental line bundle of the minimal resolution after base change:
\begin{equation}
    \label{eq:fundlinebundleafterbc}
12 \mathrm{deg}\; \Ll_{\tilde{S}/C} = \nu_0'(\tilde{\Delta})+\nu_\infty'(\tilde{\Delta}) + p \ep
\end{equation}

We refer the reader to Table 5.2 of \cite{shiodaschuttmordellweillatticesbook} for a complete description of fiber type changes after a ramified base change - we will be able to avoid doing many easy, but tedious, computations using information from that table.

To check whether a pair $(S,\gamma)$ is $\Ll$-stable, we have to determine whether the quantity (\ref{eq:fundlinebundleafterbc}) is equal to $12d$.

{\lem{\label{lem:epsilonissmall} Suppose $(S,\gamma)$ is $\Ll$-stable.
Then:
\begin{equation}
    \label{eq:epsmall}
    \ep \leq \frac{12 d}{p}
\end{equation}
}}

\begin{proof}
After base changing, the fibers which do not lie over the ramification locus contribute $p\ep$ to the degree of the minimal discriminant. 

If $(S,\gamma)$ is $\Ll$-stable, then the degree of the minimal discriminant is $12d$, which means:
\[ p\ep \leq 12d \quad \implies \quad \ep \leq \frac{12d}{p} \]
\end{proof}

The upper bound on $\ep$ gives us a lower bound for $\nu_0(\Delta)+\nu_\infty(\Delta)$ that depends on $d$ and $p$:
\[ \ep = 12 d - (\nu_0(\Delta)+\nu_\infty(\Delta)) \leq \frac{12d}{p} \quad \implies \quad 12d \left(1 - \frac{1}{p} \right) \leq \nu_0(\Delta)+\nu_\infty(\Delta) \]

We can obtain a bound on $\nu_0(\Delta)+\nu_\infty(\Delta)$ that \emph{only} depends on $d$ by taking $p = 5$ in the previous equation:
\begin{equation}
    \label{eq:discat0inftybig}
    \nu_0(\Delta)+\nu_\infty(\Delta) \geq \frac{48d}{5}
\end{equation}

{\prp{Suppose $(S,\gamma)$ is $\Ll$-stable for some $p \geq 5$. 
Then $d \leq 2$.}}

\begin{proof}

By the pigeonhole principle, we must have either $\nu_0(\Delta) \geq \frac{24d}{5}$ or $\nu_\infty(\Delta) \geq \frac{24d}{5}$.
WLOG, we may assume that $\nu_0(\Delta) \geq \frac{24d}{5}$.

Now, if $d \geq 3$, that means $\nu_0(\Delta) \geq \frac{72}{5} > 14$,
so we must have a fiber of type $I_m$ or $I_n^*$ over $0$.
Furthermore:
\begin{itemize}
    \item In the $I_m$ case, we have $m \geq \frac{24d}{5}$.
    \item In the $I_n^*$ case, we have $n \geq \frac{24d-30}{5}$
\end{itemize}

After base changing, an $I_m$ fiber becomes a fiber of type $I_{pm}$.
Since $pm \geq p \cdot \frac{24d}{5} \geq 24d > 12d$, this situation never leads to an $\Ll$-stable pair.

Similarly, if we have a fiber of type $I_n^*$, then after base changing, it will contribute:
\[ 6+ p\cdot \frac{24d-30}{5} \geq 6 +(24d-30) =24(d-1) = 12 d + 12 d - 24 = 12d + (12d-24) \]
Thus, the new $I_n^*$ fiber's contribution to the degree of the discriminant will be strictly greater than $12d$ whenever $d> 2$.
Since the discriminant has degree exactly $12d$ for $\Ll$-stable pairs,
this shows there are no $\Ll$-stable pairs with $d>2$.

\end{proof}

It is now easy to classify fiber configurations that can appear on an $\Ll$-stable rational elliptic surface ($d =1$)
or a K3 surface $(d = 2$). 

\subsection{Rational surfaces}

Let $S\to \PP^1$ be a rational elliptic surface,
and $\ep = 12 - (\nu_0(\Delta)+\nu_\infty(\Delta))$.
We will classify fiber configurations that can appear on $\Ll$-stable pairs for each possible value of $\ep$, and show that:
\begin{itemize}
    \item There are $\Ll$-stable pairs $(S,\gamma)$ with $\ep = 0$.
    The index of $\gamma$ can be arbitrarily large in this case.
    \item The other possibilities for $\ep$ are 1,2. 
    We will give examples of $\Ll$-stable pairs with $\ep = 1$ which have index 5/7.
    We will give examples where $\ep = 2$ which are of degree 5. 
\end{itemize}

\subsubsection{$\ep = 0$}

We start by classifying pairs $(S,\gamma)$ with $\ep = 0$.
We will show that there are only 4 possible fiber configurations:
\begin{itemize}
    \item If $\ep = 0$, then $S\to \PP^1$ has at most 2 singular fibers.
    \item We will show that $S\to \PP^1$ must have \emph{at least} 2 singular fibers if $S$ is rational,
    and if we have exactly 2 fibers, they must both have additive reduction.
    \item Finally, we will show that a rational elliptic surface with a fiber of type $I_n^*$,
    with $n\neq 0$, can never be part of an $\Ll$-stable pair of index $p\geq 5$.
    
\end{itemize}

It will then follow from the classification of singular fibers that \emph{if} $(S,\gamma)$ is $\Ll$-stable,
$S$ is rational, $\ep = 0$ and $\gamma$ has degree $p\geq 5$,
then the fiber configruatin of $S$ is one of the following:
\[ I_0^*+I_0^* \quad IV^*+IV \quad III^*+III \quad II^*+II \]

Furthermore:
\begin{itemize}
    \item These fibers configurations can all be realized on a rational elliptic surface. 
    \item If we take any Galois cover of index $p$ branched over the two points with singular fibers,
    then after base changing and minimizing, we end up with an isomorphic rational elliptic surface.
\end{itemize}

{\lem{Let $S\to \PP^1$ be a rational elliptic surface,
and assume $(S,\gamma)$ is $\Ll$-stable of prime index $p\geq 5$.
Then $S\to \PP^1$ does not have any fibers of type $I_n^*$ with $n \neq 0$.

}}
\begin{proof}
Suppose we have a fiber of type $I_n^*$ which does not lie over the ramification locus of $\gamma$.
Then after base changing, we have $p$ fibers of type $I_n^*$.
But a rational elliptic surface can only have one fiber of type $I_n^*$ with $n \neq 0$.

If we have an $I_n^*$ fiber over the ramification locus, then after base change,
it becomes an $I_{pn}^*$ fiber.

If $n = 1, p = 5$, we obtain an $I_5^*$, which has 10 components and contributes a degree 11 factor to the discriminant.
This can't happen on a rational elliptic surface, since by Shioda-Tate, that would imply the Picard number of the surface is at least $2 + (10-1) = 11>10$.

For any other combination of $n, p$, we end up with a fiber that has even more components, so we could never have $\Ll$-stability.

\end{proof}

{\lem{Let $S\to \PP^1$ be a rational elliptic surface.
Then $S$ has at least 2 singular fibers.}}

\begin{proof}

we know that the degree of $\Delta$ is 12, and the only fibers that contribute 12 to the degree of the discriminant are fibers of type $I_{12}$ and fibers of type $I_6^*$.
\begin{itemize}
    \item If we had a fiber of type $I_{12}$, then Shioda-Tate would say:
    \[\rho(S) = 2+\rank (E(K)) + (12-1)  \neq 10 \]
    \item Fibers of type $I_{6}^*$ have $11$ components, so again, Shioda-Tate would imply $\rho(S)>10$.
\end{itemize}
\end{proof}

{\lem{Let $S\to \PP^1$ be a rational elliptic surface,
and assume that $S\to \PP^1$ has exactly 2 singular fibers.
Then both fibers have additive reduction.}}

\begin{proof}

We may assume the fibers are at 0,$\infty$.
Let $m = \nu_0(\Delta)$ and $n = \nu_\infty(\Delta)$.
Then:
\begin{itemize}
    \item The degree of the discriminant is $12$, which gives us the constraint $m+n = 12$.
    \item The number of components in the fiber over $0$ is $m-1$ if we have additive reduction over 0 and $m$ if we have a fiber of type $I_m$.
    Similarly, we have $n-1$ components in the fiber over $\infty$ when the fiber has additive reduction and $n$ components if we have a fiber of type $I_n$.
\end{itemize}

But Shioda-Tate tells us:
\[ 10 = rank (E(K))+\left (m-\frac{1}{2} \pm \frac{1}{2}\right)+\left (n-\frac{1}{2} \pm \frac{1}{2}\right) \]
This is only possible if the number of components over $0$ is $m-1$ and the number of components over $\infty$ is $n-1$,
i.e. we have additive reduction over the two fibers.

\end{proof}

{\prp{Let $S\to \PP^1$ be a rational elliptic surface with $\ep = 0$.
Then the fiber configuration of $S\to \PP^1$ is one of:
\[ II^*+II, \quad III^*+III ,\quad IV^*+IV ,\quad I_0^*+I_0^* .\]
}}

\begin{proof}

First, since $\ep = 0$, we have no singular fibers away from $0,\infty$.
Thus, $S\to \PP^1$ has no more than 2 singular fibers.
Since $S$ is rational, that means we have exactly 2 singular fibers and they both have additive reduction.
Furthermore, since $(S,\gamma)$ is $\Ll$-stable, we can't have a fiber of type $I_n^*$ with $n \neq 0$.

Thus, the fibers over 0,$\infty$ have to be chosen from:
\[ II ,\quad III ,\quad IV, \quad I_0^* ,\quad IV^* ,\quad III^* ,\quad II^*. \]
Since, $\nu_0(\Delta) +\nu_\infty(\Delta) = 12$ and the fiber type determines $\nu_0(\Delta), \nu_\infty (\Delta)$, an easy computation shows that the only combinations that work are the ones listed above.

\end{proof}

\subsubsection{$\ep \neq 0$}

If $(S,\gamma)$ is $\Ll$-stable of index $p\geq 5$ and $\ep \neq 0$,
then we must have $p\ep \leq 12$,
so either $\ep = 2, p = 5$ or $\ep = 1$ and $p \leq 11$.

That means:
\[\nu_0(\Delta) + \nu_\infty(\Delta) \geq 10 \]
Furthermore,
after base changing and minimizing,
the contribution from the fibers over 0 and $\infty$ must be \emph{strictly less} than the contribution before base changing:
since $\ep \neq 0$, there are fibers away from the ramification locus whose contribution is \emph{guaranteed} to increase after base changing.

If we have a fiber of type $I_m$ or $I_n^*$ over the ramification locus,
then their contribution to the discriminant can only increase after base change.
Similarly, if we have a fiber of type $II/III/IV$, then a base change of prime index $p\geq 5$ will either leave the isomorphism type of the fiber unchanged, or else fibers of type $II$ can become fibers of type $II^*$/fibers of type $III$ can become fibers of type $III^*$/fibers of type $IV$ can become fibers of type $IV^*$. 

Thus, if $\ep \neq 0$ and $(S,\gamma)$ is $\Ll$-stable,
then one of the points in the ramification locus has a fiber of type $IV^*/III^*/II^*$.

We need two more lemmas:

{\lem{Let $(S,\gamma)$ be $\Ll$-stable, rational of index $p \geq 5$.
If $S$ has exactly one additive fiber,
then that fiber is of type $II^*$ and $p = 5$.
}}

\begin{proof}

Since $S$ is rational, if $S$ has exactly one additive fiber,
then $S$ has at least 3 singular fibers,
so in particular, there must be fibers of type $I_m$ somewhere.

Furthermore, since $(S,\gamma)$ is $\Ll$-stable, $S\to \PP^1$ does not have fibers of type $I_n^*$ ($n\geq 1$),
so the fiber configuration consists of a single fiber chosen from:
\[ I_0^*, \quad  II,\quad III, \quad IV, \quad IV^*, \quad III^*, \quad II^*. \]

Now, let $\ep_{ss}$ be the total contribution to the degree of the discriminant from the semistable fibers.
Then $\ep_{ss} \geq 2$, and is equal to 2 if and only if the additive fiber is of type $II^*$.
Furthermore, if $\ep_{ss}>2$,
then after base change, the semistable fibers contribute $p\ep_{ss} \geq 15$ to the degree of the discriminant.
But that would mean $(S,\gamma)$ is not $\Ll$-stable.

Thus, $\Ll$-stability either forces $S \to \PP^1$ to have two additive fibers,
or a type $II^*$ fiber.

\end{proof}

{\lem{Let $S\to \PP^1$ be a rational surface with a fiber of type $IV^*$ and a fiber of type $II$,
and let $\gamma : \PP^1\to \PP^1$ be a Galois cover of index $p\geq 5$.
Then the minimal resolution of the base change is not rational.
}}

\begin{proof}

First, we may assume that fiber of type $IV^*$ must lie over the ramification locus of $\gamma$ - otherwise,
after base change, we would have $p$ fibers of type $IV^*$,
so there is no chance that the surface is rational.

Now, we consider two cases:
\begin{itemize}
    \item If the type $II$ fiber lies over the second ramification point,
    then after base change, the fiber configuration over $0,\infty$ is either $II^*+IV$ (if $p\equiv 5\pmod 6$) or $II+IV^*$ (if $p \equiv 1 \pmod 6$).
    
    If we have $II+IV^*$, then those fibers alone contribute $14$ to the degree of the discriminant, so we don't have a rational surface. Otherwise, the configuration is of type $IV^*+II$, so those two fibers still contribute 10 to the degree of $\Delta$.
    
    But the discriminant also has a factor of degree $2p$ after base changing, since there were singular fibers which did not lie over the discriminant locus.
    Thus, we see that the surface after base change can not be rational if we have fibers of type $IV^*,II$ over the ramification locus.
    
    \item If the type $II$ fiber does not lie over the ramification point, then it splits into $p$ distinct fibers of type $II$. Since $p \geq 5$, they contribute a factor of degree at least 10 to the degree of the discriminant after base change.
    
    But we still have a fiber of type $IV$ after base change, which will contribute an extra factor of degree 4 to the discriminant, so altogether the discriminant has degree at least 14.
\end{itemize}

Thus, in every possible case, we find that the surface after base change is not rational if the fiber configuration contains $IV^*+II$ and $p\geq 5$.

\end{proof}

\begin{itemize}
    \item If we have a fiber of type $II^*$, it must be paired with a fiber of type $II$ or a pair of $I_1$s.
    (Recall that $II^*+I_2$ is not possible).
    
    \item If we have a fiber of type $III^*$, it must be paired with another additive fiber.
    Thus, we either have $III^*+III$ or $III^*+II+I_1$.
    
    \item If we have a fiber of type $IV^*$, it must be paired with another additive fiber.
    Furthermore, the second additive fiber can't be of type $II$.
    Thus, we either have $IV^*+IV$ or $IV^*+III+I_1$.
\end{itemize}

\begin{table}[]
\begin{tabular}{|l|l|l|l|l|}
\hline
No. & $\mathrm{deg} (\gamma)$ & $\pi^{-1}(0)$ & $\pi^{-1}(0)$ & Remaining fibers \\ \hline
1   & Any $p \geq 5$           & $IV^*$        & $IV$          &                  \\ \hline
2   & Any $p \geq 5$           & $III^*$       & $III$         &                  \\ \hline
3   & Any $p \geq 5$           & $II^*$        & $II$          &                  \\ \hline
4   & 5                        & $II^*$        & $I_0$         & $II$             \\ \hline
5   & 5                        & $II^*$        & $I_0$         & $I_1+I_1$         \\ \hline
6   & 5                        & $II^*$        & $I_1$         & $I_1$        \\ \hline
7   & 5                        & $IV^*$        & $III$         & $I_1$            \\ \hline
8   & 7                        & $III^*$       & $II$          & $I_1$            \\ \hline
\end{tabular}
\vspace{4mm}
\caption{Fiber configurations that occur on $\Ll$-stable pairs with $S$ rational}
\label{tab:configs}
\end{table}

In the next section, we will analyze the rational surfaces with these configurations in more detail to determine whether the associated fibrations are rank/$WC$-stable .

\subsection{K3 surfaces}

We briefly discuss $\Ll$-stability for K3 surfaces,
since there are fewer pairs when $S\to \PP^1$ is a K3 surface.
{\prp{Let $(S,\gamma)$ be an $\Ll$-stable pair of prime index $p\geq 5$ and assume $\Ll_{S/\PP^1} \cong \Oo(2)$.
Then:
\begin{itemize}
    \item The index is exactly 5.
    \item $S\to \PP^1$ has fibers of type $II^*$ over the two points in the ramification locus of $\gamma$.
\end{itemize}
}}

\begin{proof}

By Lemma \ref{lem:epsilonissmall}, we have $\ep \leq \frac{24}{p} \leq \frac{24}{5}<5$,
so $\ep \leq 4$ and $\nu_0(\Delta)+\nu_\infty(\Delta) \geq 20$.

By the pigeonhole principle,
we must have $\nu_0(\Delta) \geq 10$ or $\nu_\infty(\Delta) \geq 10$.

Suppose we have a fiber of type $I_m$, with $m \geq 10$,
or a fiber of type $I_n^*$, with $n \geq 4$, over 0.

In the first case, after base change, we have a fiber of type $I_{pm}$, with $pm \geq 50$,
so the degree of the discriminant is clearly much larger than 24.

If we have a fiber of type $I_n^*$ with $n \geq 4$, then after base change,
we have a fiber of type $I_{pn}^*$, with $pn \geq 20$.
Since fibers of type $I_{pn}^*$ contribute $6 + pn$ to the degree of the discriminant and $6 + 20 > 24$,
we deduce:
\begin{itemize}
    \item The fiber over 0 must be of type $II^*$, because that is the only fiber type that contributes 10 to the degree of the discriminant, but is not of type $I_m$ or $I_n^*$.
    \item Thus:
    \[ \nu_0(\Delta)+\nu_\infty(\Delta) = 10 + \nu_\infty(\Delta)  \geq 20 \quad \implies \quad \nu_\infty (\Delta) \geq 10 \]
    \item We now repeat the argument we use to rule out a fiber of type $I_m$ or $I_n^*$ over 0 to show that the fiber over $\infty$ is also of type $II^*$.
\end{itemize}

\end{proof}

There are not that many possible configurations one can obtain on a K3 with two fibers of type $II^*$.
However, they are all $\Ll$-stable with respect to the cover $\gamma: \PP^1\to \PP^1$ of degree 5 which ramifies over the points with a type $II^*$ fiber:
the type $II^*$ fibers become fibers of type $II$ after base changing,
so we have $\nu_0'(\tilde{\Delta})=\nu_\infty'(\tilde{\Delta}) =2$ and $p\ep = 20$.

\section{Rank and $WC$-stability}
\label{sec:examples}

When $S\to \PP^1$ is a rational elliptic surface,
we can determine the Mordell-Weil group from the fiber configuration.
Furthermore, we can use Theorem 8.33 from \cite{shiodaschuttmordellweillatticesbook} to obtain explicit generators for the Mordell-Weil group.
Thus, we can actually determine which of the $\Ll$-stable configurations we found in the previous section are actually rank stable/$WC$-stable. 

\subsection{Extremal to Extremal}
\label{ssec:ext2ext}

First, we show that nothing interesting happens when we have an $\Ll$-stable pair which appears in rows 1,2, or 3 in Table \ref{tab:configs}.
These are the pairs where we have a fiber configuration of type $IV^*+IV/III^*+III/II^*+II$ over the two ramification points of $\gamma$.

After base changing and minimizing, we end up with a rational elliptic surface with the same fiber configuration.\footnote{The location of the fibers might change after base change, e.g. it's possible to start with $III$ over 0 and $III^*$ over $\infty$, and end up with an elliptic surface with a type $III^*$ over 0 and a type $III$ over $\infty$.}
Since these configurations force the elliptic surface to be extremal (see, e.g. \cite{mirandaperssonextremalrational}) and the fiber configurations don't change after base change, it's clear that all $\Ll$-stable pairs of this form are rank stable.

\subsection{Extremal to non-Extremal}
\label{ssec:ext2nonext}

Next, we discuss the $\Ll$-stable pairs numbered $4/5/6$ in Table \ref{tab:configs}.
\begin{itemize}
    \item These fiber configurations force the original elliptic surface to be extremal (see \cite{mirandaperssonextremalrational}).
    \item We will show that we can base change by a degree 5 cover to obtain a new rational elliptic surface which is not extremal.
\end{itemize}

Thus, we can apply Theorem \ref{thm:ext2nonext} to compute the kernel of the Weil-Chatelet groups.

Let $\gamma: \PP^1 \to \PP^1$ be a Galois cover of degree $p\geq 5$ ramified over $0,\infty$,
and let $S\to \PP^1$ be a rational elliptic surface with a fiber of type $II^*$ over 0.

\begin{itemize}
    \item If we have exactly one other singular fiber (necessarily of type II) and that fiber lies over $\infty$,
    then as we just saw in (\ref{ssec:ext2ext}), the pair is $\Ll$-stable and rank stable.
    
    \item If we have exactly one other singular fiber which does \emph{not} lie over the ramification locus,
    then after base change, we either have a fiber of type $II^*$ and $p$ fibers of type $II$, or we have $p+1$ fibers of type II.
    In order to obtain a rational elliptic surface, we would need $p+1 = 6$, i.e. $p = 5$.
    Since a fiber of type $II^*$ becomes a fiber of type $II$ when we pass to a degree 5 extension,
    everything works out.
    Furthermore, $II^*+II$ is extremal, and all fibers after base change are irreducible,
    so the rank goes from 0 to 8 after base change.

    \item Next, we assume we have fibers $II^*+2I_1$ to start.
    After base changing, we will either end up with $II+10 I_1$ (if the two $I_1$ fibers do not lie over the ramification locus) or $II+5I_1 + I_5$ (if one of the two $I_1$ fibers lies over the ramification locus).
    The rank after base change is 8 in the $II+10I_1$ case and the rank is 4 in the $II+5I_1 + I_5$ case.

\end{itemize}

Note that none of these elliptic surfaces have torsion after base change,
so we can apply Theorem \ref{thm:ext2nonext} to deduce that these pairs are not $WC$-stable.
In fact, we know that $H^1(G,E_\gamma(K)) \cong \ZZ/5\ZZ$ if the rank after base change is 4, and $\ZZ/5\ZZ \times \ZZ/5 \ZZ$ if the rank is 8. 

\subsection{Non-Extremal Examples}

There are two fiber configurations left to analyze.
\begin{itemize}
    \item The elliptic surfaces have rank 1 before base changing.
    \item After base changing by a degree $p$ cover (where $p = 5$ or $7$), we obtain a rational elliptic surface with Mordell-Weil group of rank $p$.
\end{itemize}

For both of these configurations, we will fix an explicit equation for the elliptic surface, and describe a generating set for the Mordell-Weil group after base change.
We will show that these pairs are $WC$-stable by applying Theorem \ref{thm:nonextremalwcstable}.
We will need explicit equations for the elliptic surfaces so that we can describe points $Q,P$ that satisfy the conditions of that theorem. 

In both of the remaining configurations, we have two additive fibers and an $I_1$ fiber.
We can always do a change of variable that places the additive fibers over $0,\infty$,
which means we can always work with an equation of the form:
\[ y^2 = x^3 + at_0^{m_0} t_1^{m_1}+bt_0^{n_0} t_1^{n_1} \qquad (a,b \in k) \]

The fiber configuration determines the exponents in the equation above, e.g. if we want fiber configuration $IV^*+III+I_1$,
we are forced to set\footnote{Of course, we can interchange the roles of $t_0, t_1$.}:
\[ y^2 = x^3 + at_0^3 t_1+bt_0^4 t_1^2 \qquad (a,b \in k) \]
Similarly, to obtain $III^*+II+I_1$, we have work with:
\[ y^2 = x^3 + at_0 t_1^3+bt_0 t_1^5 \qquad (a,b \in k) \]

We still have to choose coefficients $a,b \in k$, but the choice of $a,b$ doesn't matter if $k$ is algebraically closed.
Once we fix a choice of $a,b$, we can find the equation of the elliptic surface after base change,
and then apply (Theorem 8.33) from (\cite{shiodaschuttmordellweillatticesbook}) to obtain a generating set for the base-changed elliptic surface.

\begin{itemize}
    \item We give an equation for a rational elliptic surface with the desired fiber configuration.
    \item We use Mathematica to determine how many points are in the generating set.
    \item We describe one of those points explicitly, and then show how to obtain the remaining generators as $\ZZ$-linear combinations of the point we gave. 

\end{itemize}

This will be enough to verify that the conditions of Theorem \ref{thm:nonextremalwcstable} are satisfied.

\subsubsection{$IV^*+III+I_1$}

Next, we discuss elliptic surfaces with fiber configuration $IV^*+III+I_1$.
For concreteness, we will assume the surface is given by the equation:

\begin{equation}
    \label{eq:ell34prebc}
    E: \quad y^2 = x^3 -  t_0^3 t_1 x + t_0^4 t_1^2 
\end{equation}

which has a fiber of type $IV^*$ over $t_0 =0$ and a fiber of type $III$ over $t_1 = 0$.
Let $\gamma : \PP^1 \to \PP^1$ be a degree 5 cover ramified over $t_0 = 0$ and $t_1 =0$.

After base changing by $\gamma$, the equation becomes:
\begin{equation}
    \label{eq:ell34bc}
    E_\gamma: \quad y^2 = x^3 -  t_0^{15} t_1^5 x + t_0^{20} t_1^{10} 
\end{equation}
and after minimizing, we obtain:
\begin{equation}
    \label{eq:ell34bc}
    E_\gamma': \quad y^2 = x^3 -  t_0^3 t_1 x + t_0^2 t_1^4 
\end{equation}
Finally, we pass to the chart $t_1 =1$:
\begin{equation}
    \label{eq:ell34bcaff}
    E_\gamma': \quad y^2 = x^3 -  t^3 x + t^2 
\end{equation}

\begin{itemize}
    \item Since the elliptic surface is rational after base change, Theorem 8.33 from \cite{shiodaschuttmordellweillatticesbook} says that the Mordell-Weil group after base change is generated by the set of points that have the form $(at^2+bt+c,dt^3+et^2+ft +g)$.

    \item The Galois group acts on this set of points.
    To describe this action explicitly, we pick a generator $\si$ of the Galois group, as well as a primitive fifth root of unity $\zeta \in \CC$, and we define:
    \begin{equation}
        \label{eq:galoisaction34}
        \si \cdot (at^2+bt+c,dt^3+et^2+ft +g) = (a\zeta^3t^2+b\zeta^2 t+c\zeta ,d\zeta^2 t^3+e \zeta t^2+ft +g\zeta^4) 
    \end{equation}

    Rather than listing all 92 points, we can just give a single element from each orbit and use (\ref{eq:galoisaction34}) to obtain the other points.

    \item Using Mathematica, we find that there are 92 such points, and they can all be described explicity. However, some of the coefficients are ``complicated" algebraic integers that we would rather not write out. Instead, we will show that all 92 points can be obtained as $\ZZ$-linear combinations of Galois conjugates of a single point.

\end{itemize}

Let:
\begin{equation}
    Q_1 = (-2^{2/5} t,t-\sqrt[5]{2} t^2)
\end{equation}

Note that $Q_1$ is not fixed by $G$, so the Galois orbit of $Q_1$ contains 5 of the generators.
Furthermore, $-Q_1$ (and the Galois conjugates of $-Q_1$) have the same monomial structure as $Q$,
so they account for another 5 elements from the generating set.

A computation shows that:
\[ Tr(Q_1) = (0,t) \]
Call this point $P_0$.
Then $P_0$ and $-P_0$ account for another pair of points from the generating set.

That leaves 80 generators unaccounted for.
Since each generator gives rise to 10 elements when we negate and take Galois conjugates,
we only need to describe 8 more points in the generating set.

Define:

\begin{align*}
    Q_0 &= Q_1 - \si \cdot Q_1\\
    Q_2& = Q_1 + \si Q_1 \\
    Q_3 &= Q_2 + \si^4 Q_1 \\
    Q_4 &= (\si+\si^2+\si^3+\si^4) \cdot Q_1 = Tr(Q_1) - Q_1 = P_0 - Q_1\\
    Q_5& = Q_0 + \si Q_0 \\
    Q_6 &= Q_2 + Q_4 \\
    Q_7 &=\si^4 \cdot Q_0 +Q_2  \\
    Q_8 &= \si^3 \cdot Q_0 + Q_3
\end{align*}
Then:
\begin{itemize}
    \item Each of these points has the form $(at^2+bt+c,dt^3+et^2+ft +g) $.
    \item None of these points is fixed by $\si$, so they give us the remaining 80 points in the generating set.
\end{itemize}

Thus, we can apply Theorem \ref{thm:nonextremalwcstable} with $Q = Q_1$ to prove that this pair is $WC$-stable.

\subsubsection{$III^*+II+I_1$}
Next, we study elliptic surfaces with configuration $III^*+II+I_1$.
Note that this is the only configuration that appears on an $\Ll$-stable but not rank stable pair of index $p\geq 7$.
This elliptic surface is discussed in Theorem 3.1 of \cite{fastenbergcycliccovers1}.

We will work with the Weierstrass equation:
\begin{equation}
    \label{eq:ell7prebc}
   E: \quad y^2 = x^3 - t_0 t_1^3 x + t_0 t_1^5
\end{equation}

After base changing (but before passing to a minimal integral equation),
the equation becomes:
\begin{equation}
    \label{eq:ell7nonmin}
  E_\gamma :\quad  y^2 = x^3 -  t_0^7 t_1^{21} x + t_0^7 t_1^{35}
\end{equation}

We can use the usual change of variable:
\[ ((t_0 t_1^5)^2y)^2 = ((t_0 t_1^5)^2 x)^3- t_0^7 t_1^{21}  ((t_0 t_1^5)^2 x) + t_0^7 t_1^{35}\]
and divide through by $t_0^6 t_1^{30}$ to obtain:
\begin{equation}
\label{eq:ell7postbcmin}
 E_\gamma': \quad   y^2 = x^3 -  t_0^3 t_1x+ t_0 t_1^5
\end{equation}

The $\ZZ[G]$-module structure on $E_\gamma(K)$ is induced by the action $u\mapsto \zeta_7 u$,
where $u = \frac{t_0}{t_1}$, $\zeta_7 = e^{2\pi i/7}$,
and $E(K)$ is isomorphic to the subgroup $E_\gamma(K)^G \subset E_\gamma(K)$.

To compute the Mordell-Weil group, we work with the equation (\ref{eq:ell7postbcmin}) on the chart $t_1 \neq 0$.
We write $t = \frac{t_0}{t_1}$ for the local parameter.

Our new equation is:
\begin{equation}
\label{eq:ell7postbcmin}
 E_\gamma': \quad   y^2 = x^3 - t^3x +t
\end{equation}

Now, $E_\gamma(k(t))$ has rank 7, so we can use the stronger version of the generating theorem:
every rational elliptic surfaces with Mordell-Weil rank equal to 7 contains exactly 56 points of the form:
\[(at+b, ct^2+dt+e) \]
and these points generate the full Mordell-Weil group (see Prop. 7.12.ii in \cite{shiodaschuttmordellweillatticesbook}).

The action of the Galois group $G = \ip{\si}$ on one of these points is:
\[ \si \cdot (at + b, ct^2 + dt + e) = (\zeta^3 at +\zeta^2 b, \zeta^5 ct^2+\zeta^4 dt + \zeta^3 e) \]
where $\zeta$ is a primitive 7th root of unity.

Let $a_1 \in k$ be a root of the polynomial:
\[ a^{21} + 5a^{14} + 6 a^7 + 1 =0\]
Define:
\begin{align*}
    b_1 &=.\frac{-4845 a_1^{24}-18546717 a_1^{17}-423971443 a_1^{10}-410387643 a_1^3}{4958072} \\
    c_1 &= i\sqrt{a_1}  \\
    d_1 &= \frac{i \sqrt{a_1} \left(-4845 a_1^{23}-18546717 a_1^{16}-423971443 a_1^9-415345715 a_1^2\right)}{9916144}  \\
    e_1 &=\frac{i \sqrt{a_1} \left(-312930 a_1^{25}-1197895591 a_1^{18}-27367766068 a_1^{11}-26135476903 a_1^4\right)}{34706504}
\end{align*}

Let $Q_1 = (a_1 t+ b_1, c_1t^2 + d_1 t + e_1)$, and define:
\begin{align*}
    Q_2 &= \si^2 \cdot Q_1 - \si^4 \cdot Q_1 + \si^6 \cdot Q_1 \\
    Q_3 &=  Q_2 - \si \cdot Q_1 + \si^3 \cdot Q_1\\
    Q_4 &= Q_1 + Q_2 + Q_3
\end{align*}
These 4 points, together with their Galois conjugates, their negatives, and the negatives of their Galois conjugates give us the 56 points in the generating set.
The trace of $Q_1$ is $(u^{-2}, u^{-3})$ with respect to these coordinates.
This corresponds to the point $(1,1)$ on $E(K)$,
so we can apply Theorem \ref{thm:nonextremalwcstable} with $Q = Q_1$ to prove $WC$-stability for this pair.

The results of this section are summarized in Table \ref{tab:summary}.

\begin{table}[]
\begin{tabular}{|l|l|l|l|l|l|l|}
\hline
$\mathrm{deg} (\gamma)$ & $\pi^{-1}(0)$ & $\pi^{-1}(0)$ & Remaining fibers & $\mathrm{rank}E(K)$ & $\mathrm{rank}E(K')$ & $\ker (WC\to WC)$                                                                                       \\ \hline
Any $p \geq 5$           & $IV^*$        & $IV$          &                  & 0                   & 0                    & 0                                                                                                       \\ \hline
Any $p \geq 5$           & $III^*$       & $III$         &                  & 0                   & 0                    & 0                                                                                                       \\ \hline
Any $p \geq 5$           & $II^*$        & $II$          &                  & 0                   & 0                    & 0                                                                                                       \\ \hline
5                        & $II^*$        & $I_0$         & $II$             & 0                   & 8                    & \begin{tabular}[c]{@{}l@{}}$(\mathbb{Z}/5\mathbb{Z})^2 $\end{tabular}    \\ \hline
5                        & $II^*$        & $I_0$         & $I_1+I_1$        & 0                   & 8                    & \begin{tabular}[c]{@{}l@{}}$(\mathbb{Z}/5\mathbb{Z})^2$\end{tabular} \\ \hline
5                        & $II^*$        & $I_1$         & $I_1$            & 0                   & 4                    & \begin{tabular}[c]{@{}l@{}}$\mathbb{Z}/5\mathbb{Z}$\end{tabular}                                     \\ \hline
5                        & $IV^*$        & $III$         & $I_1$            & 1                   & 5                    & 0                                                                                                       \\ \hline
7                        & $III^*$       & $II$          & $I_1$            & 1                   & 7                    & 0                                                                                                       \\ \hline
\end{tabular}
\vspace{4mm}
\caption{Mordell-Weil ranks and Weil-Chatelet kernels for $\Ll$-stable pairs}
\label{tab:summary}
\end{table}

\section{$WC$ kernels for $\Ll$-unstable pairs over $\PP^1$}

To obtain more interesting examples, we need to drop the $\Ll$-stability condition.
While this condition \emph{is} relevant to applications,
it is very restrictive.
Furthermore, the other types of stability may be of interest to people who are not interested in $\Ll$-stability.

Now, without $\Ll$-stability, things become much more difficult:
\begin{itemize}
    \item We have not yet discussed what happens with semistable fibrations, since they are never $\Ll$-stable. Ruling out new torsion points requires a deeper argument: we can no longer appeal to the presence of an additive fiber to rule out all torsion of order $>4$.
    However, we can use modular curves to rule out new torsion points when $\gamma$ is Galois of index $p>3$.
    
    \item The real problem with pairs which are not $\Ll$-stable is understanding the structure of the Mordell-Weil group after base change: not only are we missing an analog of Theorem 8.33 in \cite{shiodaschuttmordellweillatticesbook}, which enabled us to compute explicit generators of Mordell-Weil group \emph{after} base change, but we can't even bound the rank of the Mordell-Weil group, since we can't systematically compute the Picard number of the elliptic surface after base change.

\end{itemize}

However, in some special cases\footnote{Essentially when we have some way to bound the rank of the Mordell-Weil group after base change.}, we can prove rank/$WC$-stability for $p>>0$.

\subsection{Semistable Surfaces}

Let $S\to \PP^1$ be a semistable rational elliptic surface.
Note that $S$ has at least 4 singular fibers.

Let $n_0$ be the order of vanishing of $\Delta$ at 0 and
$n_\infty$ the order of vanishing at $\infty$.
Let $n_1,n_2, \hdots, n_\ell$ be the order of vanishing of the discriminant over the remaining points in the discriminant locus.
Then:
\[ n_0 + n_\infty +\sum_{i = 1}^\ell n_i = 12 \qquad n_0, n_\infty \geq 0 \qquad n_1, n_2, \hdots, n_\ell > 0\]

Set $\ep = 2$ if $n_0,n_\infty$ are both positive, $\ep = 1$ if exactly one of $n_0, n_\infty$ is positive and $\ep = 0$ if $n_0 = n_\infty =0$.

\begin{itemize}
    \item The original elliptic surface has $\ell + \ep$ singular fibers.
    After base changing, we have $p\ell + \ep$ singular fibers.
    \item The rank of the trivial lattice before base changing is:
    \[ rank T = \sum m_\nu -1 = 12 - (\ell + \ep) \]
    After base changing, the rank of the trivial lattice is:
    \[ rank T_\gamma = 12p - (p\ell + \ep) \] 
    \item Thus, we can compute the difference between $rank T_\gamma$ and $rank T$:
    \begin{align*}
     rank T_\gamma - rank T &= (12p - (p\ell + \ep)) - ( 12 - (\ell + \ep))\\
     &=12(p-1) -(p-1)\ell - (p-1) \ep \\
     &=(p-1)(12-\ell - \ep)
    \end{align*}
\end{itemize}

Now, Shioda-Tate tells us:
\begin{align*}
    rank\; E(K) + rank T = \rho(S) = 10\\
    rank\; E(K_\gamma) + rank T_\gamma = \rho(S_\gamma) \leq 10p
\end{align*}
Subtracting the second equation from the first:
\begin{align*}
    (rank\; E(K_\gamma) - rank E(K)) - (rank T_\gamma- rank T_\gamma) = \rho(S_\gamma) -\rho(S) \leq 10p-10
\end{align*}
Rearranging, and using the previous computation, we obtain:
\begin{align*}
     (rank\; E(K_\gamma) - rank E(K)) \leq 10(p-1) - (p-1)(12-\ell - \ep)
\end{align*}
Finally, we simplify the RHS to obtain:
\begin{equation}
\label{eq:rankboundshiodatatess}
     (rank\; E(K_\gamma) - rank E(K)) \leq (p-1)(-2+\ell+\ep)
\end{equation}

Now, while we have a bound on  $(rank\; E(K_\gamma) - rank E(K))$,
the bound depends on $p$.
This means that we can't prove rank stability by arguing that the jump in rank would exceed the bound on rank from Shioda-Tate.

We can, however, bound the dimension of the kernel $WC(E/K)\to WC(E/K_\gamma)$ as an $\FF_p$ vector space,
at least if the Mordell-Weil group after base changing is torsion-free.

Assume $E(K_\gamma)$ is torsion-free.
Let $\si$ be a generator of the Galois group, and let $\Oo = \ZZ[e^{2\pi i/p}]$. 
    Then $E(K_\gamma) = E(K) \oplus M$, for some regular $\ZZ[G]$-module $M$.
     Let $\tau = \sum_{\si^\ell \in G} \si^\ell \in \ZZ[G]$, and $M_\tau \subset M$ be the submodule of $M$ annihilated by $\tau$.
    Then $M_\tau$ is an $\Oo$-module.
    
     Let $M_0 = (1-\si)M$. Note that $M_0 \cap M^G = M_0 \cap E(K) = 0$, so the map $m\mapsto (1-\si) m$ is injective. Thus, $rank M_0 = rank M$ as an abelian group.
    Furthermore, $M_0 \subset M_\tau$, so $M_0$ is also a regular $\Oo$-module. The kernel of $WC(E/K) \to WC(E/K_\gamma)$ is isomorphic to $M_\tau/M_0$.
    In order to bound the $\FF_p$ rank of the quotient, we need to bound the rank of $M_\tau$ as an $\Oo$-module.

    The rank of $M_\tau$ as an $\Oo$-module is $(p-1)$ times the rank of $M$ as an abelian group. Furthermore, the rank of $M$ as an abelian group is $rank(E(K_\gamma)) - rank(E(K))$.
    Thus, we can apply (\ref{eq:rankboundshiodatatess}) to deduce that the rank of $M_\tau$ as an $\Oo$-module is bounded above by $-2+\ell +\ep$, with equality only if $\rho(S_\gamma) = 12p = h^{1,1}(S_\gamma)$.

\subsection{Delsartes and Fastenberg Surfaces}

To understand $WC$-kernels in general, we need a better understanding of Mordell-Weil groups of general elliptic surfaces.
At the very least, we need some way of bounding the rank of the Mordell-Weil group \emph{after} base changing.
\begin{itemize}
    \item For rational elliptic surfaces, we really just needed to know the fiber configuration to determine the exact rank of the Mordell-Weil group.

    \item Mordell-Weil groups of K3 surfaces are not as predictable as those of rational elliptic surfaces,
    but they are much better understood than Mordell-Weil groups of general elliptic surfaces. We have a uniform bound on the rank, but the Mordell-Weil rank is no longer determined by the fiber configuration alone, and we have no analog of the Generator Theorem (8.33 in \cite{shiodaschuttmordellweillatticesbook}).

\end{itemize}

For elliptic surfaces of Kodaira dimension 1,
    we only have comparable results for special families. 

A Delsartes elliptic surface over $k$ is an elliptic surface over $\PP^1$ whose generic fiber
is given by an equation of the form:
\[ y^2 = x^3 + at^m x + bt^n \]
where $a,b \in k$ and $m,n \in \NN$.
If $S\to \PP^1$ is a Delsartes surface and $\gamma: \PP^1\to \PP^1$ is a Galois cover branched over $0,\infty$, then the base $S_\gamma \to \PP^1$ is also a Delsartes surface.

    Furthermore, if $E/k(t)$ is the generic fiber of a Delsartes surface, then rank $E(k(t))$ is bounded above by 68.
    The bound,is sharp, e.g. the elliptic surface:
    \[ y^2 = x^3 + t^{360} x+ 1 \]
    has rank 68.
Altogether, this gives us an upper bound on the \emph{change} in rank after we base change. 
On the other hand,
we know that the change in rank is a multiple of $p-1$.
Thus, we have rank stability if we base change by a Galois extension of index $p > 68$.
See (\cite{heijnedelsarte},\cite{shiodaexplicitalgorithm}) for more on these surfaces.

Similarly, the families of elliptic surfaces studied in \cite{fastenbergcycliccovers1}, \cite{fastenbergmordellweil} 
are stable under base change.
Furthermore, for each of those families, there is a bound on the rank of the Mordell-Weil group which is independent of the degree of the base extension. Thus, we can argue as we did with the Delsartes surfaces to obtain rank stability for cyclic extensions of sufficiently large prime degree.

\end{document}